\documentclass[reqno]{amsart}
\usepackage{amssymb}

\usepackage[all]{xy}
\CompileMatrices
\newdir{ >}{{}*!/-10pt/\dir{>}}

\theoremstyle{plain}

\theoremstyle{definition}

\newtheorem{emptythm}{}[section]

\newcommand{\too}{\longrightarrow}

\newcommand{\ul}[1]{\underline{#1}}

\DeclareMathOperator{\Hom}{{\rm Hom}}
\DeclareMathOperator{\MotHom}{{\rm hom}}
\DeclareMathOperator{\MotEnd}{{\rm end}}
\DeclareMathOperator{\HM}{{\rm H}}

\DeclareMathOperator{\res}{{\rm res}}

\DeclareMathOperator{\Ker}{{\rm Ker}}

\DeclareMathOperator{\id}{{\rm id}}
\DeclareMathOperator{\Gal}{{\rm Gal}}

\DeclareMathOperator{\khar}{{\rm char}}
\DeclareMathOperator{\Spec}{Spec}

\DeclareMathOperator{\Pic}{{\rm Pic}}

\DeclareMathOperator{\Br}{{\rm Br}}

\DeclareMathOperator{\CHTor}{{\rm A}}
\DeclareMathOperator{\et}{{\rm et}}

\DeclareMathOperator{\Chow}{\mathfrak{Chow}}         
\DeclareMathOperator{\CH}{{\rm CH}}                  

\newcommand{\Ch}[3]{\CH_{#1}(#2)_{#3}}                
\newcommand{\Mhom}[4]{\MotHom_{#1}(#2,#3)_{#4}}       
\newcommand{\Mend}[3]{\MotEnd_{#1}(#2)_{#3}}          

\newcommand{\MK}{{\rm K}^{M}}                        
\newcommand{\QK}{{\rm K}}                            
\newcommand{\CM}{{\rm M}}                            

\newcommand{\Tm}[1]{\ul{#1}}

\newcommand{\Gm}{\mathbb{G}_{{\rm m}}}

\newcommand{\Z}{\mathbb{Z}}
\newcommand{\N}{\mathbb{N}}
\renewcommand{\P}{\mathbb{P}}

\newcommand{\cf}{\textsl{cf.}\ }
\newcommand{\eg}{\textsl{e.g.}\ }
\newcommand{\ie}{\textsl{i.e.}\ }

\begin{document}

\title[Permutation modules and Chow motives]
{Permutation modules and Chow motives of geometrically rational surfaces}

\author{Stefan Gille}
\address{Stefan Gille, Department of
         Mathematical and Statistical Sciences,
         University of Alberta
         Edmonton T6G 2G1, Canada}

\email{gille@ualberta.ca}

\thanks{This work has been supported by
        an NSERC research grant}


\date{May 15, 2014}

\begin{abstract}
We prove that the Chow motive with integral coefficient of a geometrically
rational surfaces~$S$ over a perfect field~$k$ is zero dimensional if and
only if the Picard group of~$\bar{k}\times_{k}S$, where~$\bar{k}$
is an algebraic closure of~$k$, is a direct summand of a
$\Gal (\bar{k}/k)$-permutation module, and~$S$ possesses a zero cycle
of degree one. As shown by Colliot-Th\'el\`ene in a letter to the author
(which we have reproduced in the appendix) this is in turn equivalent
to~$S$ having a zero cycle of degree~$1$ and $\CH_{0}(k(S)\times_{k}S)$
being torsion free.
\end{abstract}

\maketitle

\section{Introduction}
\label{IntroSect}\bigbreak

\noindent
Let~$k$ be a perfect field with algebraic closure~$\bar{k}$,
and~$S$ a geometrically rational $k$-surface. Set $S:=\bar{k}\times_{k}S$.
The Picard group~$\Pic\bar{S}$ is
then a free $\Z$-module of finite rank and also a discrete
$G_{k}:=\Gal (\bar{k}/k)$-module. It is called a $G_{k}$-permutation
module if it has a $\Z$-basis which is permuted by the
action of~$G_{k}$, and $G_{k}$-invertible if it is a direct summand
of a $G_{k}$-permutation module.

\medbreak

The aim of this short note is to prove the following result.

\medbreak

\noindent
{\bf Theorem A.}
{\it
Let~$S$ be a geometrically rational $k$-surface. Then the motive of~$S$
in the category of effective Chow motives $\Chow (k)$ with integral coefficients
is $0$-dimensional if and only if~$\Pic\bar{S}$ is an invertible $G_{k}$-module
and~$S$ has a $0$-cycle of degree~$1$.
}

\medbreak

\noindent
Recall that a motive~$M$ in~$\Chow (k)$ is called $0$-dimensional
if it is a direct summand of a finite direct sum of Tate twists
of motives of $k$-\'etale algebras.

\smallbreak

Our proof of this theorem is constructive, and can be used to get
motivic decompositions of motives of such surfaces. We illustrate
this in the (rather easy) case of Del-Pezzo surfaces of degree~$5$.
As another byproduct of our proof we get the (well known?) fact
that a motive in the category of Chow motives with rational coefficients
is zero dimensional if and only if it is geometrically split.

\smallbreak

Being a $0$-dimensional motive is stable under field extensions, and
so if~$S$ is a geometrically rational surface whose motive is $0$-dimensional
in~$\Chow (k)$ then the abelian group~$\CH_{0}(l\times_{k}S)$ is torsion free
for all field extensions~$l\supseteq k$. Colliot-Th\'el\`ene has shown
in a letter~\cite{BriefCT13} to the author, which we have reproduced
in the Appendix, that if~$S$ is a geometrically rational surface with
$0$-cycle of degree~$1$ and no torsion in~$\CH_{0}(k(S)\times_{k}S)$, where~$k(S)$
denotes the function field of~$S$, then~$\Pic\bar{S}$ is a invertible
$G_{k}$-module. This together with Theorem~A proves the converse of the
assertion above.

\smallbreak

In Sections~\ref{CH0-SubSect} and~\ref{universaltorsorCor}
we extend these equivalences further:

\medbreak

\noindent
{\bfseries Theorem B.}
{\it
Let~$S$ be a geometrically rational surface over the
perfect field~$k$, and~$T^{S}$ the $k$-torus with character
group~$\Pic\bar{S}$. Then the following four assertions are equivalent:
\begin{itemize}
\item[(i)]
The motive of~$S$ in~$\Chow (k)$ is $0$-dimensional.

\smallbreak

\item[(ii)]
The $G_{k}$-module~$\Pic\bar{S}$ is invertible and~$S$ has a $0$-cycle
of degree~$1$.

\smallbreak

\item[(iii)]
The group $\CH_{0}(k(S)\times_{k}S)$ is torsion free
and~$S$ has a zero cycle of degree~$1$. (Note that
for a geometrically rational surface~$S$ the group~$\CH_{0}(S)$
is torsion free if and only if the degree map $\CH_{0}(S)\too\Z$
is injective.)

\smallbreak

\item[(iv)]
The degree map~$\CH_{0}(l\times_{k}Y^{c})\too\Z$ is an isomorphism for all
smooth compactifications~$Y^{c}$ of~$T^{S}$-torsors~$Y$ over~$S$ and all
field extensions~$l$ of~$k$.
\end{itemize}
}

\medbreak

\noindent
One may wonder here, whether there exists an example of a geometrically rational
$k$-surface~$S$, whose motive is zero dimensional, but which has no $k$-rational
point.

\bigbreak

\noindent
{\bfseries Acknowledgement.}
I would like to thank Professor Colliot-Th\'el\`ene for allowing me to reproduced his
letter~\cite{BriefCT13} as an appendix, for carefully answering a lot of my
questions about geometrically rational surfaces, and also for corrections and
useful comments.

\smallbreak

Further I would like to thank Jochen Kuttler, J\'an Min\'a\v{c}, Nikita Semenov, and Sasha Vishik
for many useful discussions.

\bigbreak

\bigbreak\goodbreak
\section{Preliminaries and notations}
\label{PreliminariesSect}\bigbreak

\begin{emptythm}
\label{PerNotationsSubSect}
{\it Notations and conventions.}
Throughout this work~$k$ denotes a perfect field with algebraic (and so
also separable) closure~$\bar{k}$. We denote the absolute Galois group
of the field~$k$ by~$G_{k}$. If~$M$ is a continuous $G_{k}$-module
we denote by~$\HM^{i}(k,M)$, or by~$\HM^{i}(G_{k},M)$, the Galois
cohomology of~$M$.

\smallbreak

If~$l/k$ is a field extension and~$X$ a $k$-scheme
we set $X_{l}:=l\times_{k}X$, and $\bar{X}:=X_{\bar{k}}$.
Similarly we set~$A_{l}:=l\otimes_{k}A$
and~$\bar{A}:=A_{\bar{k}}$ if~$A$ is a commutative $k$-algebra.

\smallbreak

If~$X$ is an integral $k$-scheme we denote the function field
by~$k(X)$.
\end{emptythm}

\begin{emptythm}
\label{ChowGroupSubSect}
{\it Chow groups.}
Let~$K$ be a field and ~$X$ a $K$-scheme of finite type. We denote by~$\CH_{i}(X)$
the Chow group of dimension $i$-cycles modulo rational equivalence as
defined in Fulton's book~\cite{IntTheory}. We denote the class of an
$i$-dimensional subvariety~$Z$ in~$\CH_{i}(X)$ by~$[Z]$.

\smallbreak

If~$R$ is a commutative ring
(with unit~$1$) we set $\Ch{i}{X}{R}:=R\otimes_{\Z}\CH_{i}(X)$,
\ie~$\Ch{i}{X}{R}$ is the Chow group of dimension $i$-cycles with
coefficients in~$R$. In case~$X=\Spec A$ is an affine scheme we also use the notation
$\Ch{i}{A}{R}$ instead of~$\Ch{i}{X}{R}$.

\smallbreak

For later use we state the following well known lemma, which \eg follows
from Karpenko and Merkurjev~\cite[Cor.\ RC.12]{KaMe13}.

\smallbreak

\noindent
{\bfseries Lemma.}
{\it
Let~$K$ be an arbitrary field, and~$X$ and~$Y$ smooth and projective
$K$-schemes which are $K$-birational to each other. Then $\deg :\CH_{0}(X)\too\Z$
is an isomorphisms if and only if $\deg :\CH_{0}(Y)\too\Z$ is an isomorphism.
}
\end{emptythm}

\begin{emptythm}
\label{ComDiagSubSect}
{\it A commutative diagram.}
Let~$X$ be a $k$-scheme and~$\tau\in G_{k}$. Then~$\tau$ induces a morphism
of schemes $\tilde{\tau}:\Spec\bar{k}\too\Spec\bar{k}$. which in turn
induces by base change the automorphism $\tilde{\tau}\times\id_{X}$ of
$\bar{X}=\Spec\bar{k}\times_{k}X$. The pull-backs of these morphisms define
then the action of~$G_{k}$ on~$\Ch{i}{\bar{X}}{R}$ by the rule
$\tau .\alpha:=(\tilde{\tau}\times\id_{X})^{\ast}(\alpha)$.
Note that $(\tilde{\tau}\times\id_{X})^{\ast}=(\tilde{\tau}^{-1}\times\id_{X})_{\ast}$.
As $\Ch{i}{\bar{X}}{R}$ is the direct limit of all $\Ch{i}{X_{l}}{R}$,
where~$l$ runs through all finite subextensions $k\subseteq l\subseteq\bar{k}$,
this makes~$\Ch{i}{\bar{X}}{R}$ a continuous $G_{k}$-module.

\smallbreak

Assume now that~$X$ is geometrically integral and smooth over~$k$. Let
$k\subseteq l\subseteq\bar{k}$ be a finite field extension of~$k$.
We identify the absolute Galois group~$G_{l}$ of~$l$ with a
closed subgroup of~$G_{k}$, and denote by~$\Gamma$ the set of
all $k$-embeddings $l\too\bar{k}$. Every~$\gamma\in\Gamma$
induces a morphism of $k$-schemes
$\tilde{\gamma}:\bar{X}\too\bar{k}\times_{k}X_{l}=\bar{k}\times_{k}l\times_{k}X$, and
$$
\bigcup\limits_{\gamma\in\Gamma}\;\tilde{\gamma}\,\; :\;\bigcup\limits_{\gamma\in\Gamma}\bar{X}
\,\too\,\bar{k}\times_{k}X_{l}
$$
is an isomorphism of $k$-schemes, where~$\bigcup\limits_{\gamma\in\Gamma}\bar{X}$
is the disjoint union of $|\Gamma|=[l:k]$-copies of~$\bar{X}$. The induced isomorphism
$$
(\tilde{\gamma}^{\ast})_{\gamma\in\Gamma}\, :\;
\Ch{i}{\bar{k}\times_{k} X_{l}}{R}\,\xrightarrow{\;\simeq\;}\,
\bigoplus\limits_{\gamma\in\Gamma}\Ch{i}{\bar{X}}{R}
$$
is $G_{k}$-linear for the following (by this isomorphism induced)
$G_{k}$-action on the direct sum $\bigoplus\limits_{\gamma\in\Gamma}\Ch{i}{\bar{X}}{R}$: For
$\tau\in G_{k}$ the $\gamma_{0}$-component of $\tau.(\alpha_{\gamma})_{\gamma\in\Gamma}$
is given by $\tau.\alpha_{\tau^{-1}\gamma_{0}}$.

\medbreak

Let~$\rho_{\gamma}\in G_{k}$ be an extension of~$\gamma\in\Gamma$ to~$\bar{k}$.
We have then a commutative diagram
\begin{equation}
\label{CH-DiagEq}
\xymatrix{
\Ch{i}{X_{l}}{R} \ar[d]_-{\res_{\bar{k}/l}} \ar[rr]^-{\res_{\bar{k}/k}}
      & & \Ch{i}{\bar{k}\times_{k}X_{l}}{R}^{G_{k}}
     \ar[d]^-{(\tilde{\gamma}^{\ast})_{\gamma\in\Gamma}}_-{\simeq}
\\
\Ch{i}{\bar{X}}{R}^{G_{l}} \ar[rr]_-{\simeq} & &
{\Big[\bigoplus\limits_{\gamma\in\Gamma}\Ch{i}{\bar{X}}{R}\Big]^{G_{k}}} \rlap{\, ,}
}
\end{equation}
where the bottom arrow maps $\alpha\in\Ch{i}{\bar{X}}{R}^{G_{l}}$
to $(\rho_{\gamma}.\alpha)_{\gamma\in\Gamma}$. Note that this map
does not depend on the choice of the extensions~$\rho_{\gamma}$.
As indicated in the diagram the bottom arrow as well as the arrow
on the right hand side are isomorphisms.
\end{emptythm}

\begin{emptythm}
\label{ChowMotSubSect}
{\it Chow motives.}
Let~$X,Y$ be smooth and projective $k$-schemes.
A correspondence of degree~$0$ with coefficients in~$R$
from~$X$ to~$Y$ is an element~$\alpha$ of
$\bigoplus\limits_{i=1}^{r}\Ch{\dim X_{i}}{X_{i}\times_{k}Y}{R}$, where
$X_{1},\ldots ,X_{r}$ are the irreducible components
of~$X$. We indicate this by $\alpha:X\vdash Y$. The composition of
two correspondences is denoted by~$\beta\circ\alpha$.

\smallbreak

Smooth projective varieties with correspondences of degree~$0$ (with
coefficients in~$R$) as morphisms and disjoint unions as direct sums
are an additive category, the category of correspondences of degree~$0$
with~$R$ as coefficient ring. The idempotent completion of this category
is the category of (effective) Chow motives over~$k$ with coefficient ring~$R$.
We denote this category~$\Chow (k,R)$ or $\Chow (k)$ if~$R=\Z$. The objects
of~$\Chow (k,R)$ are pairs $(X,\pi)$, where~$X$ is a smooth projective
$k$-scheme and $\pi: X\vdash X$ is an idempotent correspondence of degree~$0$,
\ie $\pi\circ\pi=\pi$. By some abuse of notation we denote the motive of
a smooth projective scheme~$X$ by the same letter~$X$.

\smallbreak

If~$M,N$ are motives, \ie objects of~$\Chow (k,R)$, we denote the group of
morphisms between them by $\Mhom{k}{M}{N}{R}$, respectively by $\MotHom_{k}(M,N)$
if~$R=\Z$. Similarly $\Mend{k}{M}{R}$ and $\MotEnd_{k}(M)$ denote the
endomorphism group of the motive~$M$ in $\Chow (k,R)$ and $\Chow (k)$,
respectively.

\smallbreak

If $l/k$ is a field extension we denote the base change functor
$\Chow (k,R)\too\Chow (l,R)$ by~$\res_{l/k}$, and set $M_{l}:=\res_{l/k}(M)$
and $\alpha_{l}:=\res_{l/k}(\alpha)$ for a motive~$M$ and a
morphism~$\alpha$ in $\Chow (k,R)$. If~$l=\bar{k}$ is the algebraic
closure we use then also~$\bar{M}$ and~$\bar{\alpha}$ for~$M_{\bar{k}}$
and~$\alpha_{\bar{k}}$, respectively.
\end{emptythm}

\begin{emptythm}
\label{SplitMotSubSect}
{\it Split motives.}
We denote the Tate motive in~$\Chow (k,R)$ by~$\Tm{R} (1)$, and set
$\Tm{R} (i):=\Tm{R} (1)^{\otimes\, i}$ and $M(i):=M\otimes\Tm{R} (i)$ for
a motive~$M$ and an integer~$i\geq 0$. In particular, $\Tm{R}=\Tm{R}(0)=\Spec k$
is the motive of the point. A motive~$M$ is called {\it split} if~$M$ is isomorphic
in~$\Chow (k,R)$ to a finite direct sum of twisted Tate motives, and {\it geometrically
split} if~$\bar{M}=M_{\bar{k}}$ is split.
Recall, see \eg Merkurjev~\cite[Prop.\ 1.5]{Me08a}, that the motive
of a smooth projective and integral scheme~$X$ is split in~$\Chow (k,R)$
if and only if
$$
\Ch{i}{X}{R}\otimes_{R}\Ch{\dim X-i}{X}{R}\,\too\, R\, ,\;
\alpha\otimes\beta\,\longmapsto\,\deg (\alpha\cap\beta)
$$
is a perfect duality for all $0\leq i\leq\dim X$. In particular, there are bases
$u_{1}^{i},\ldots ,u_{n_{i}}^{i}$ and $v_{1}^{i},\ldots ,v_{n_{i}}^{i}$
of the free $R$-modules $\Ch{i}{X}{R}$ and~$\Ch{\dim X-i}{X}{R}$,
respectively, such that
$$
\deg (u_{r}^{i}\cap v_{s}^{i})\, =\,\delta_{rs}\, ,
$$
where~$\delta_{rs}$ is the Kronecker delta. Moreover, in this case we have
$$
\id_{X}\, =\,\sum\limits_{i=0}^{\dim X}\sum\limits_{j=1}^{n_{i}}
v^{i}_{j}\times u^{i}_{j}\;\in\,\Mend{k}{X}{R}\, ,
$$
and the product of cycles induces an isomorphism
$$
\bigoplus\limits_{i=0}^{\dim X}\Ch{i}{\bar{X}}{R}\otimes_{R}\Ch{\dim X-i}{\bar{X}}{R}
\,\xrightarrow{\;\simeq\;}\,\Ch{\dim X}{\bar{X}\times_{\bar{k}}\bar{X}}{R}\,
=\,\Mend{\bar{k}}{\bar{X}}{R}\, .
$$

\smallbreak

\noindent
{\bfseries Example.}
Let~$S$ be a geometrically rational $k$-surface, \ie a smooth projective
geometrically integral $k$-scheme of dimension~$2$, such that $\bar{S}=\bar{k}\times_{k}S$
is rational. The motive of such a surface~$S$ is geometrically split
in~$\Chow (k)$ and so also in~$\Chow (k,R)$ for all coefficient rings~$R$,
see \eg~\cite[Thm.\ 2.6]{Gi10}.
\end{emptythm}

\begin{emptythm}
\label{Mot-EtSubSect}
{\it Motives of \'etale algebras and permutation modules.}
We denote by~$RG_{k}$ the group ring of~$G_{k}$ over~$R$.

\smallbreak

A continuous $RG_{k}$-module~$M$ is called
an {\it $RG_{k}$-permutation module} if~$M$ is a free and
finitely generated $R$-module which has a basis which is
permuted by the action of~$G_{k}$. If the coefficient ring
is~$\Z$ we call a $\Z G_{k}$-permutation module a {\it $G_{k}$-permutation
module}. A direct summand of a $RG_{k}$-permutation module
is called an {\it invertible $RG_{k}$-module}, respectively an {\it invertible
$G_{k}$-module} if~$R=\Z$.

\medbreak

An example of a $RG_{k}$-permutation module is~$\Ch{0}{\bar{E}}{R}$
for~$E$ a $k$-\'etale algebra, and in fact every $RG_{k}$-permutation
module is isomorphic to~$\Ch{0}{\bar{E}}{R}$ for some $k$-\'etale
algebra~$E$.

\smallbreak

If $\alpha :\Spec E\,\vdash\,\Spec F$ is a correspondence of degree~$0$
between the spectra of two $k$-\'etale algebras then
$$
\alpha_{\ast}\, :\;\Ch{0}{\bar{E}}{R}\,\too\,\Ch{0}{\bar{F}}{R}\, ,\;
\beta\,\longmapsto\,\alpha\circ\beta
$$
is a $RG_{k}$-morphism. This defines a functor from the subcategory
of~$\Chow (k,R)$ generated by \'etale $k$-algebras into the category
of $RG_{k}$-permutation modules which is full and faithful,
\cf~\cite[Sect.\ 7]{ChMe06}.
\end{emptythm}

\begin{emptythm}
\label{0-dimSubSect}
{\it $0$-dimensional motives.}
A motive~$M$ in~$\Chow (k,R)$ is called {\it $0$-dimensional} if
there exist $k$-\'etale algebras $E_{0},E_{1},\ldots ,E_{n}$, such
that~$M$ is a direct summand of the motive
$$
\bigoplus\limits_{i=0}^{n}(\Spec E_{i})(i)\, .
$$
Note that if the motive of a smooth projective scheme~$X$ is
$0$-dimensional in $\Chow (k)=\Chow (k,\Z)$ then it is also
$0$-dimensional in~$\Chow (k,R)$ for all coefficient rings~$R$.

\medbreak

\noindent
{\bfseries Example.}
The motive of a $k$-rational surface~$S$ is $0$-dimensional in~$\Chow (k,R)$
for every (commutative) coefficient ring~$R$. In fact,
since~$k$ is perfect there exists a smooth projective surface~$S'$, and
morphisms $S'\too S$ and $S'\too\P^{2}_{k}$ which are compositions of blow-ups
in closed points, see~\cite[Chap.\ II, Cor.\ 21.4.2]{Cubic}. Therefore by
the  blow-up formula (\cf~\cite{Ma68} or~\cite[Chap.\ 6]{IntTheory}) the motive
of~$S$ is a direct summand of the motive of~$S'$ which in turn by the blow-up
formula and the projective bundle theorem is isomorphic to
$$
\Tm{R}\oplus\,\bigoplus\limits_{i=1}^{m}(\Spec l_{i})(1)\,\oplus\Tm{R}(2)
$$
for some finite field extensions $l_{i}\supseteq k$.

\bigbreak

For later use we mention the following easy lemma.

\smallbreak

\noindent
{\bfseries Lemma.}
{\it
Let~$X$ be a smooth projective geometrically integral $k$-scheme
whose motive is $0$-dimensional in~$\Chow (k,R)$. Then
$\deg:\Ch{0}{X}{R}\too R$ is surjective. In particular, if moreover~$R$ is torsion
free and the motive of~$X$ in~$\Chow (k,R)$ is geometrically split then
the degree map is an isomorphism.
}

\begin{proof}
Since~$X$ is geometrically integral the field~$k$ is algebraically
closed in its function field~$k(X)$, and therefore
$$
\res_{k(X)/k}\, :\;\Mhom{k}{F}{L}{R}\,\too\,
\Mhom{k(X)}{k(X)\otimes_{k}F}{k(X)\otimes_{k}L}{R}
$$
is an isomorphism for all
finite field extensions $F,L$ of~$k$. This implies since the motive of~$X$
in~$\Chow (k,R)$ is a direct summand of motives of Tate twists of separable
field extensions of~$k$ that the restriction map
$\Ch{0}{X}{R}\too\Ch{0}{X_{k(X)}}{R}$ is an isomorphism. The first claim follows.

\smallbreak

The last assertion is a consequence of the fact that the kernel of
the degree map is torsion if the motive of~$X$ in $\Chow (k,R)$ is
geometrically split.
\end{proof}
\end{emptythm}

\begin{emptythm}
\label{Z-trivialSubSect}
{\it $\Z$-trivial varieties.}
There is a closely related property of schemes. A geometrically integral $k$-scheme~$X$ is
called {\it $\Z$-trivial} if $\deg:\CH_{0}(X_{l})\too\Z$ is an isomorphism for all field
extensions $l\supseteq k$ with $X(l)\not=\emptyset$, and
{\it strongly $\Z$-trivial} if $\deg:\CH_{0}(X_{l})\too\Z$ is an isomorphism
for all field extensions $l\supseteq k$.

\smallbreak

\noindent
{\bfseries Remark.}
The notion $\Z$-trivial, or more general $R$-trivial for an
arbitrary coefficient ring~$R$, has been introduced by Karpenko
and Merkurjev~\cite[Def.\ 2.3]{KaMe13}.

\medbreak

As shown by Merkurjev~\cite[Thm.\ 2.11]{Me08b}
a smooth projective $k$-scheme~$X$ which is geometrically integral is strongly $\Z$-trivial
if and only if the class of the generic point in~$\CH_{0}(X_{k(X)})$ is defined
over~$k$, and this is in turn equivalent to the assertion that for any cycle
module~$\CM_{\ast}$ in the sense of Rost~\cite{Ro96} the cycle cohomology
group~$\HM^{0}(X,\CM_{\ast})$ is equal~$\CM_{\ast}(k)$ via the canonical
map $\CM_{\ast}(k)\too\HM^{0}(X,\CM_{\ast})\subseteq\CM_{\ast}(k(X))$.

\medbreak

\noindent
{\bfseries Example.}
Let~$X$ be a geometrically rational smooth and projective $k$-scheme whose
motive is $0$-dimensional in~$\Chow (k)$. Then by the lemma
in~\ref{0-dimSubSect} above the degree map~$\CH_{0}(X_{l})\too\Z$ is an isomorphism
for all $l\supseteq k$, \ie~$X$ is strongly $\Z$-trivial.

\medbreak

This property has the following consequence.
\end{emptythm}

\begin{emptythm}
\label{HS-SpSeqLem}
{\bfseries Lemma.}
{\it
Let~$X$ be a geometrically integral smooth and projective $k$-scheme,
which is strongly $\Z$-trivial. Then $\Br (k)=\HM^{2}(k,\bar{k}^{\times})\too\HM^{2}_{\et}(X,\Gm)$
is injective, and so the natural homomorphism
$\Pic X\too (\Pic\bar{X})^{G_{k}}$ is an isomorphism.
}

\begin{proof}
The latter statement follows from the first by the Hochschild-Serre
spectral sequence.

\smallbreak

For the first assertion it is enough to show that the canonical map
$\Br (k)\too\Br (k(X))$ is injective. Since~$k$ is by assumption perfect
there are no elements of exponent $p=\khar k >0$ in~$\Br (k)$ if~$k$ has
positive characteristic~$p$ by a theorem of Albert~\cite[p.\ 109]{StAlg}.
Hence it is enough to show that for~$n$ prime
to~$\khar k$ the natural homomorphism $\HM^{2}(k,\mu_{n})\too\HM^{2}(k(X),\mu_{n})$
is injective, where $\mu_{n}\subseteq\bar{k}^{\times}$ denotes the group
of $n$th roots of unity. For such~$n$ the assignment
$l\mapsto\HM^{\ast}(l,\mu_{n}^{\otimes\,\ast -1})$ is a cycle module
in the sense of Rost, see~\cite{Ro96}, and so by the above mentioned result
of Merkurjev~\cite[Thm.\ 2.11]{Me08b} the natural homomorphism
$\HM^{2}(k,\mu_{n})\too\HM^{2}(k(X),\mu_{n})$ is injective.
\end{proof}
\end{emptythm}

\bigbreak\goodbreak
\section{Maps between \'etale algebras and smooth projective schemes in
the category of Chow motives}
\label{MapsEt-SchSect}\bigbreak

\begin{emptythm}
\label{Et-SchSubSect}
{\it Morphisms from motives of Tate twists of \'etale algebras to motives of smooth
projective schemes.}
We fix throughout this section a coefficient ring~$R$.
Let~$X$ be a smooth projective and geometrically integral $k$-scheme and~$E$
a $k$-\'etale algebra, \ie~$E=l_{1}\times\ldots\times l_{d}$ for some finite separable
extensions $l_{i}\supseteq k$.

\smallbreak

Let $0\leq i\leq n:=\dim X$ be an integer.
We have a homomorphism
$$
\psi_{E,X}^{i}\, :\;\Mhom{k}{(\Spec E)(i)}{X}{R}\,\too\,
\Hom_{RG_{k}}(\Ch{0}{\bar{E}}{R},\Ch{i}{\bar{X}}{R})\, ,\;
\alpha\,\longmapsto\,\bar{\alpha}_{\ast}\, ,
$$
where $\bar{\alpha}_{\ast}:\Ch{0}{\bar{E}}{R}\too\Ch{i}{\bar{X}}{R}$ maps
the cycle~$\beta$ to $\bar{\alpha}\circ\beta$. Note that since~$\alpha$ is
defined over~$k$ the homomorphism~$\bar{\alpha}_{\ast}$ is $RG_{k}$-linear.

\smallbreak

Let $\bar{k}\supseteq l\supseteq k$ be an algebraic extension of~$k$. The absolute
Galois group~$G_{l}$ of~$l$ can be identified with a closed subgroup of~$G_{k}$
and so any~$RG_{k}$-homomorphism is also a $RG_{l}$-homomorphism. The following
diagram commutes:
\begin{equation}
\label{Et-SchEq}
\xymatrix{
\Mhom{l}{(\Spec E_{l})(i)}{X_{l}}{R} \ar[rr]^-{\psi_{E_{l},X_{l}}^{i}} & &
\Hom_{RG_{l}}(\Ch{0}{\bar{E}}{R},\Ch{i}{\bar{X}}{R})
\\
\Mhom{k}{(\Spec E)(i)}{X}{R} \ar[u]^-{\res_{l/k}} \ar[rr]^-{\psi_{E,X}^{i}} & &
\Hom_{RG_{k}}(\Ch{0}{\bar{E}}{R},\Ch{i}{\bar{X}}{R}) \ar[u]_-{\subseteq} \rlap{\, .}
}
\end{equation}

\medbreak

\noindent
{\bfseries Lemma.}
{\it
Let~$X$ and~$E=l_{1}\times\ldots\times l_{d}$ be as above. Assume
that the natural homomorphism
$$
\res_{\bar{k}/l_{t}}\, :\;\Ch{i}{X_{l_{t}}}{R}\,\too\,\Ch{i}{\bar{X}}{R}^{G_{l_{t}}}
$$
is an isomorphism for all $t=1,\ldots, d$, where $G_{l_{t}}\subseteq G_{k}$
is the absolute Galois group of~$l_{t}$. Then~$\psi_{E,X}^{i}$ is an
isomorphism.
}

\begin{proof}
This is clear if~$k=\bar{k}$ is an algebraically closed field. For the general
case we can assume that $E=l$ is a finite separable field extension of~$k$.
Using the commutative diagram~(\ref{Et-SchEq}) we see then that it is enough
to show that the restriction map
$$
\res_{\bar{k}/k}\, :\;\Ch{i}{X_{l}}{R}=\Mhom{k}{l}{X}{R}\too\Mhom{\bar{k}}{\bar{l}}{\bar{X}}{R}
=\Ch{i}{l\times_{k}\bar{X}}{R}
$$
is an isomorphism onto the $G_{k}$-invariant elements in~$\Mhom{l}{\bar{l}}{\bar{X}}{R}$.
But this follows from the commutative diagram~(\ref{CH-DiagEq}) in~\ref{ComDiagSubSect} and the
assumption that $\res_{\bar{k}/l}:\Ch{i}{X_{l}}{R}\too\Ch{i}{\bar{X}}{R}^{G_{l}}$
is an isomorphism.
\end{proof}
\end{emptythm}

\begin{emptythm}
\label{Sch-EtSubSect}
{\it Morphisms from motives of smooth projective schemes to motives of Tate twists of \'etale
algebras.} We continue with the notation of the last subsection.

\smallbreak

We have a homomorphism

$\;\;\phi_{X,E}^{i}\, :\;$
$$
\qquad\Mhom{k}{X}{(\Spec E)(i)}{R}\,\too\,
\Hom_{RG_{k}}(\Ch{i}{\bar{X}}{R},\Ch{0}{\bar{E}}{R})\, ,\;
\alpha\,\longmapsto\,\bar{\alpha}_{\ast}\, .
$$
This homomorphism commutes also with base change, \ie
if~$l\supseteq k$ is a field extension then the diagram
\begin{equation}
\label{Sch-EtEq}
\xymatrix{
\Mhom{l}{X_{l}}{(\Spec E_{l})(i)}{R} \ar[rr]^-{\phi_{X_{l},E_{l}}^{i}} & &
\Hom_{RG_{l}}(\Ch{i}{\bar{X}}{R},\Ch{0}{\bar{E}}{R})
\\
\Mhom{k}{X}{(\Spec E)(i)}{R} \ar[u]^-{\res_{l/k}} \ar[rr]^-{\phi_{X,E}^{i}} & &
\Hom_{RG_{k}}(\Ch{i}{\bar{X}}{R},\Ch{0}{\bar{E}}{R}) \ar[u]_-{\subseteq} \rlap{\, .}
}
\end{equation}
commutes.

\medbreak

\noindent
{\bfseries Lemma.}
{\it
Let $E=l_{1}\times\ldots\times l_{d}$ be as above.
Assume that the motive of~$X$ in~$\Chow (k,R)$ is geometrically split
and that
$$
\res_{\bar{k}/l_{t}}\, :\;\Ch{n-i}{X_{l_{t}}}{R}\,\too\,\Ch{n-i}{\bar{X}}{R}^{G_{l_{t}}}
$$
is an isomorphism for the finite field extensions
$\bar{k}\supset l_{t}\supseteq k$ with absolute Galois group
$G_{l_{t}}\subseteq G_{k}$, $t=1,\ldots ,d$. Then~$\phi_{E,X}^{i}$
is an isomorphism.
}

\begin{proof}
If~$k=\bar{k}$ is algebraically closed this follows since the motive
of~$\bar{X}$ is split and so (\cf~\ref{SplitMotSubSect})
$$
\Ch{i}{\bar{X}}{R}\times\Ch{n-i}{\bar{X}}{R}\,\too\, R\, ,\; (\alpha,\beta)
\,\longmapsto\,\deg (\alpha\cap\beta)
$$
is a perfect pairing. As in the proof of the lemma in~\ref{Et-SchSubSect}
this implies the general case using now the commutative diagram~(\ref{Sch-EtEq})
instead of~(\ref{Et-SchEq}).
\end{proof}
\end{emptythm}

\begin{emptythm}
\label{Mot-PerModThm}
{\bfseries Theorem.}
{\it
Let~$X$ be a smooth and projective $k$-scheme of dimension~$n$,
whose motive is geometrically split in~$\Chow (k,R)$. Assume that there
is an $k$-\'etale algebra~$E$ and correspondences of degree~$0$
$$
\alpha\, :\; (\Spec E)(i)\,\vdash\, X\quad\mbox{and}\quad
\beta\, :\; X\,\vdash\, (\Spec E)(i)
$$
with corresponding $RG_{k}$-linear maps
$$
f\, =\,\psi_{E,X}^{i}(\alpha):\Ch{0}{\bar{E}}{R}\too
\Ch{i}{\bar{X}}{R}
$$
and
$$
g\, =\,\phi_{X,E}^{i}(\beta):
\Ch{i}{\bar{X}}{R}\too\Ch{0}{\bar{E}}{R}
$$
for some integer $0\leq i\leq n=\dim X$.

\smallbreak

Then we have $f\cdot g=\id_{\Ch{i}{\bar{X}}{R}}$ if and
only if
\begin{equation}
\label{a-b1Eq}
\bar{\alpha}\circ\bar{\beta}\, =\,\sum\limits_{j=1}^{m}v_{j}\times u_{j}\, ,
\end{equation}
where $u_{1},\ldots ,u_{m}\in\Ch{i}{\bar{X}}{R}$ and
$v_{1},\ldots ,v_{m}\in\Ch{n-i}{\bar{X}}{R}$ are dual basis with respect
to the perfect pairing $\Ch{i}{\bar{X}}{R}\times\Ch{n-i}{\bar{X}}{R}\too R$,
$(\gamma,\rho)\mapsto\deg (\gamma\cap\rho)$, see~\ref{SplitMotSubSect}.
}

\begin{proof}
We have a $\bar{k}$-algebra isomorphisms
$\bar{E}\xrightarrow{\simeq}\bar{k}^{r}$ with corresponding open immersions
$p_{j}:\Spec\bar{k}\too\Spec\bar{E}$. We can then write
$\bar{\alpha}=\sum\limits_{j=1}^{r}
(p_{j}\times\id_{\bar{X}})^{\ast}(\alpha_{j})$ and $\bar{\beta}=
\sum\limits_{j=1}^{r}(\id_{\bar{X}}\times p_{j})^{\ast}(\beta_{j})$
for appropriate $\alpha_{j}$'s in~$\in\Ch{i}{\bar{X}}{R}$ and
$\beta_{j}$'s in~$\Ch{n-i}{\bar{X}}{R}$, respectively.
In terms of the dual basis we have $\alpha_{j}=
\sum\limits_{s=1}^{m}a_{sj}\cdot u_{s}$ and
$\beta_{j}=\sum\limits_{t=1}^{m}b_{jt}\cdot v_{t}$
for some $a_{sj},b_{jt}\in R$. Since $\bar{\alpha}
\circ\bar{\beta}$ is in the image of the injective homomorphism
$\Ch{i}{\bar{X}}{R}\otimes_{R}\Ch{n-i}{\bar{X}}{R}\too\Ch{n}{\bar{X}\times_{\bar{k}}\bar{X}}{R}$,
$\gamma\otimes\rho\mapsto\gamma\times\rho$, and this image is a free
$R$-module with basis $u_{i}\times v_{j}$, $1\leq i,j\leq m$, by the K\"unneth
isomorphism, see~\ref{SplitMotSubSect}, the equation~(\ref{a-b1Eq})
is equivalent to
\begin{equation}
\label{a-bEq}
\sum\limits_{j=1}^{r}a_{sj}\cdot b_{jt}\, =\,\delta_{st}
\end{equation}
for all $1\leq s,t\leq m$, where~$\delta_{st}$ denotes the Kronecker delta.

\smallbreak

On the other hand the cycles~$\alpha_{j}$ and~$\beta_{j}$ correspond
by~\ref{Et-SchSubSect} and~\ref{Sch-EtSubSect} to maps
$f_{j}$~in $\Hom_{R}(\Ch{0}{\bar{k}}{R},\Ch{i}{\bar{X}}{R})$
and $g_{j}$ in $\Hom_{R}(\Ch{i}{\bar{X}}{R},\Ch{0}{\bar{k}}{R})$,
respectively, such that
$$
f=\sum\limits_{j=1}^{r}f_{j}\qquad\mbox{and}\qquad
g=\sum\limits_{j=1}^{r}g_{j}
$$
via the isomorphism $\Ch{0}{\bar{E}}{R}\simeq\Ch{0}{\bar{k}}{R}^{r}$
induced by the pull-backs along the open immersions $p_{1},\ldots ,p_{r}$.

\smallbreak

We identify now $\Ch{i}{\bar{X}}{R}$ with
$\Hom_{R}(\Ch{0}{\bar{k}}{R},\Ch{i}{\bar{X}}{R})$ naturally, and the group
$\Ch{n-i}{\bar{X}}{R}$ with $\Hom_{R}(\Ch{i}{\bar{X}}{R},\Ch{i}{\bar{k}}{R})$
via the intersection product using that the motive
of~$\bar{X}$ is split. Then we have
$f_{j}=\sum\limits_{s=1}^{m}a_{sj}\cdot u_{s}$ and
$g_{j}=\sum\limits_{t=1}^{m}b_{jt}\cdot v_{t}$, and therefore
$f(g(u_{l}))=\sum\limits_{j=1}^{r}f_{j}(g_{j}(u_{l}))=
\sum\limits_{s=1}^{m}\big(\,\sum\limits_{j=1}^{r}a_{sj}\cdot b_{jl}\,\Big)
\cdot u_{s}$
for all~$1\leq l\leq m$. Hence since $u_{1},\ldots ,u_{m}$ is an $R$-basis
of~$\Ch{i}{\bar{X}}{R}$ the equation $f\cdot g=\id_{\Ch{i}{\bar{X}}{R}}$
is equivalent to~(\ref{a-b1Eq}).
\end{proof}

\medbreak

This has the following consequence.

\smallbreak

\noindent
{\bfseries Corollary.}
{\it
Let~$X$ be a smooth and projective $k$-scheme of dimension~$n$,
whose motive is geometrically split in~$\Chow (k,R)$. Assume there are
$k$-\'etale algebras $E_{0},\ldots ,E_{n}$ and correspondences of
degree~$0$
$$
\alpha\, =\,\sum\limits_{i=0}^{n}\alpha_{i}\, :\,\bigoplus\limits_{i=0}^{n}
(\Spec E_{i})(i)\,\vdash\, X
$$
and
$$
\beta\, =\,\big(\,\beta_{i}\,\big)_{i=0}^{n}\, :\; X\,\vdash\,
\bigoplus\limits_{i=0}^{n} (\Spec E_{i})(i)
$$
with corresponding $RG_{k}$-linear maps
$f_{i}=\psi_{E_{i},X}^{i}(\alpha_{i}):\Ch{0}{\bar{E}_{i}}{R}\too
\Ch{i}{\bar{X}}{R}$ and $g_{i}=\phi_{X,E_{i}}^{i}(\beta_{i}):
\Ch{i}{\bar{X}}{R}\too\Ch{0}{\bar{E}_{i}}{R}$ for $i=0,\ldots ,n$.

\smallbreak

Then we have $\bar{\alpha}\circ\bar{\beta}=\id_{\bar{X}}$ in
$\Chow (\bar{k},R)$ if and only if $f_{i}\cdot g_{i}=\id_{\Ch{i}{\bar{X}}{R}}$
for all $i=0,\ldots ,n$.
}
\end{emptythm}

\begin{emptythm}
\label{zeroDimMotSubSect}
{\it Zero dimensional motives.}
Recall that a motive~$M$ in~$\Chow (k,R)$ is called {\it zero dimensional}
if it is a direct summand of $\bigoplus\limits_{i=0}^{n}(\Spec E_{i})(i)$ for
some $k$-\'etale algebras~$E_{i}$. If the coefficient ring~$R$
has the property that all projective $R$-modules are free this implies that~$M$
is also geometrically split. Hence one implication of the corollary above
can be formulated as follows:

\medbreak

\noindent
{\bfseries Corollary.}
{\it
Let~$X$ be a smooth projective $k$-scheme whose motive in~$\Chow (k,R)$ is
zero dimensional. Assume that
\begin{itemize}
\item[(i)]
the motive of~$X$ is geometrically split, or that

\smallbreak

\item[(ii)]
every projective $R$-module is free.
\end{itemize}
Then the $RG_{k}$-module $\Ch{i}{\bar{X}}{R}$ invertible
for all $0\leq i\leq \dim X$.
}
\end{emptythm}

\begin{emptythm}
\label{SplittingMapsSubSect}
{\it A splitting criterion.}
The converse of the corollary in~\ref{zeroDimMotSubSect} above seems to
be only true under some further assumptions on~$X$. Let for this~$X$ be a smooth
and projective $k$-schemes whose motive is geometrically split in
$\Chow (k,R)$. We assume further that
\begin{itemize}
\item[{\bf (GC)}]
If $\bar{k}\supseteq l\supseteq k$ is an algebraic field extension
then  the restriction homomorphism $\Ch{i}{X_{l}}{R}\too\Ch{i}{X_{\bar{k}}}{R}$
is an isomorphism onto the $G_{l}$-invariant elements of
$\Ch{i}{X_{\bar{k}}}{R}$ for all $0\leq i\leq n:=\dim X$ (where as
above we identify~$G_{l}$ with a subgroup of~$G_{k}$), and

\smallbreak

\item[{\bf (RN)}]
{\it Rost nilpotence} is true for~$X$ in $\Chow (k,R)$, \ie given a field
extension $l\supseteq k$, and a correspondence~$\alpha$ in~$\Mend{k}{X}{R}$,
such that $\alpha_{l}=0$ then~$\alpha$ is nilpotent: $\alpha^{\circ\, N}=0$
for some integer~$N\geq 1$.
\end{itemize}

\medbreak

\noindent
We choose for all $0\leq i\leq n$ a surjective $RG_{k}$-linear morphism
$P_{i}\xrightarrow{f_{i}}\Ch{i}{\bar{X}}{R}$ with~$P_{i}$
a $RG_{k}$-permutation module. We can assume, see~\ref{Mot-EtSubSect},
that there is a $k$-\'etale algebra~$E_{i}$, such that
$P_{i}=\Ch{0}{(E_{i})_{\bar{k}}}{R}$ for all $0\leq i\leq n$.
We have then the following result.

\medbreak

\noindent
{\bfseries Theorem.}
{\it
If the morphisms~$f_{i}$ are split, \ie if there are
$RG_{k}$-linear maps $g_{i}:\Ch{i}{X_{\bar{k}}}{R}\too\Ch{0}{(E_{i})_{\bar{k}}}{R}$,
such that $f_{i}\cdot g_{i}=\id_{\Ch{i}{\bar{X}}{R}}$ for all
$0\leq i\leq n$, then the motive of the scheme~$X$ is a direct summand
of $\bigoplus\limits_{i=0}^{n}\Spec (E_{i})(i)$ in the category $\Chow (k,R)$.
}

\begin{proof}
By assumption~{\bf (GC)} the lemmas in~\ref{Et-SchSubSect}
and~\ref{Sch-EtSubSect} imply the existence of cycles $\alpha_{i}\in\Ch{i}{E_{i}\times_{k}X}{R}=
\Mhom{k}{(\Spec E_{i})(i)}{X}{R}$ and $\beta_{i}\in
\Ch{n-i}{X\times_{k}E_{i}}{R}=\Mhom{k}{X}{(\Spec E_{i})(i)}{R}$, such that
$$
f_{i}=\psi_{E_{i},X}^{i}(\alpha_{i})\qquad\mbox{and}\qquad
g_{i}=\phi_{X,E_{i}}^{i}(\beta_{i})
$$
for $i=0,\ldots ,n=\dim X$. By the corollary to Theorem~\ref{Mot-PerModThm}
together with~\ref{SplitMotSubSect} we have then
$$
\sum\limits_{i=0}^{n}\bar{\alpha}_{i}\circ\bar{\beta}_{i}\, =\,\id_{\bar{X}}\, .
$$
Since Rost nilpotence is true for~$X$ in $\Chow (k,R)$ by assumption~{\bf (RN)}
we get therefore $\sum\limits_{i=0}^{n}\alpha_{i}\circ\beta_{i}=\id_{X}+\gamma$,
where~$\gamma$ is a nilpotent correspondence of degree~$0$ on~$X$. Hence the claim.
\end{proof}
\end{emptythm}

\begin{emptythm}
\label{ratCoeffExpl}
{\bfseries Example.}
Let~$R=F$ be a field of characteristic~$0$. Then for any field
extension $l/k$ and any finite type $k$-scheme~$Y$ the base
change homomorphism $\Ch{i}{Y}{F}\too\Ch{i}{Y_{l}}{F}$ is
injective for all $i\in\N$ and so Rost nilpotence is true for
all motives in $\Chow (k,F)$. Moreover, if $\bar{k}\supseteq l\supseteq k$
is an algebraic field extension and~$Y$ is geometrically integral
then the base change map $\Ch{i}{Y_{l}}{F}\too\Ch{i}{\bar{Y}}{F}^{G_{l}}$
is an isomorphism for all $0\leq i\leq\dim Y$, where $G_{l}\subseteq G_{k}$
denotes the absolute Galois group of~$l$.

\smallbreak

Since by Maschke's theorem every surjective $G_{k}$-linear map
between continuous and $F$-finite dimensional $G_{k}$-modules splits
we conclude from the corollary in~\ref{zeroDimMotSubSect} and the
theorem in~\ref{SplittingMapsSubSect} the following (well known?) fact.

\medbreak

\noindent
{\bfseries Corollary.}
{\it
Let~$F$ be a field of characteristic~$0$, and~$X$ a smooth and
projective $k$-scheme. Then the motive of~$X$ in~$\Chow (k,F)$
is zero dimensional if and only if it is geometrically split.
}
\end{emptythm}

\bigbreak\goodbreak
\section{Geometrically rational surfaces with zero dimensional Chow motive}
\label{GeomRatSurfSect}\bigbreak

\begin{emptythm}
\label{GeomRatSurfSubSect}
{\it Geometrically rational surfaces.}
Let~$S$ be a geometrically rational $k$-surface.
Recall that by~\cite{Gi10,Gi14} Rost nilpotence is true
for~$S$ in~$\Chow (k)$.

\medbreak

\noindent
{\bfseries Theorem.}
{\it
Let~$S$ be a geometrically rational $k$-surface.
Then the motive of~$S$ in~$\Chow (k)$ is zero dimensional
if and only if~$\Pic\bar{S}$ is an invertible $G_{k}$-module
and~$S$ has a zero cycle of degree~$1$.
}

\begin{proof}
One direction is a consequence of the corollary to Theorem~\ref{Mot-PerModThm}
and the lemma in~\ref{0-dimSubSect}.

\smallbreak

To prove the converse, we show first that if the motive of~$S$ is $0$-dimensional
then~$S$ is strongly $\Z$-trivial. By Merkurjev~\cite[Thm.\ 2.11]{Me08b},
see~\ref{Z-trivialSubSect}, it is enough to show that the class of the generic
point in~$\CH_{0}(S_{k(S)})$ is defined over~$k$. Since~$S$ has by assumption a
$0$-cycle of degree one this follows if we show that~$\CH_{0}(S_{k(S)})$ is torsion free.

\smallbreak

Let~$T^{S}$ be the $k$-torus with character group~$\Pic\bar{S}$,
\ie $T^{S}=\Spec (\bar{k}[\Pic\bar{S}]^{G_{k}})$, and set $K:=\bar{k}(S)$.
Since~$S_{\bar{k}(S)}$ is rational we have by~\cite[Thm.\ 0.1]{Bl81}
an exact sequence (identifying $\Gal (K/k(S))=G_{k}$)
$$
\HM^{1}\big(G_{k},\MK_{2}(K(S))/\MK_{2}(k(S))\big)\,\too\,
\CHTor_{0}(S_{k(S)})\,\too\,\HM^{1}(G_{k},T^{S}_{k(S)})\, ,
$$
where $\CHTor_{0}(S_{k(S)})=\Ker\big(\deg:\CH_{0}(S_{k(S)})\too\Z\big)$
is the torsion part of~$\CH_{0}(S_{k(S)})$.

\smallbreak

Since~$\Pic\bar{S}$ is $G_{k}$-invertible the torus~$T^{S}$ is a direct
factor of a quasi-trivial torus and thus by Hilbert's Theorem~90 we
have~$\HM^{1}(G_{k},T^{S}_{k(S)})=0$. The group on the left hand side vanishes
by the main result of Colliot-Th\'el\`ene~\cite{CT83}. Hence~$\CHTor_{0}(S_{k(S)})=0$,
and therefore~$S$ is strongly $\Z$-trivial.

\medbreak

It follows then that $\res_{\bar{k}/l}:\CH_{0}(S_{l})\too
\CH_{0}(\bar{S})=\CH_{0}(\bar{S})^{G_{l}}$ is an isomorphism,
and by Lemma~\ref{HS-SpSeqLem} that $\Pic S_{l}\too (\Pic\bar{S})^{G_{l}}$
is an isomorphism for all intermediate fields extensions
$k\subseteq l\subseteq\bar{k}$, where $G_{l}\subseteq G_{k}$ denotes the
absolute Galois group of~$l$. The same is obviously true
for~$\CH_{2}$, \ie $\CH_{2}(S_{l})\xrightarrow{\simeq}\CH_{2}(\bar{S})^{G_{l}}$,
and so we can apply the theorem in~\ref{SplittingMapsSubSect} which
shows that the motive of~$S$ in~$\Chow (k)$ is $0$-dimensional as claimed.
We are done.
\end{proof}
\end{emptythm}

\begin{emptythm}
\label{ComputationSubSect}
{\it A computational remark.}
Let~$S$ be a geometrically rational $k$-surface. Assume that the
motive of~$S$ in~$\Chow (k)$ is $0$-dimensional. Then by the
lemma in~\ref{0-dimSubSect} the $G_{k}$-module~$\CH_{0}(\bar{S})$
is isomorphic to the trivial $G_{k}$-module~$\Z$, and the
same is true for~$\CH_{2}(\bar{S})$. Since the trivial $G_{k}$-module~$\Z$
corresponds to the motive~$\Tm{\Z}$ this implies by the corollary to
Theorem~\ref{Mot-PerModThm} and the theorem in~\ref{SplittingMapsSubSect}
that the motive of~$S$ is a direct summand of
\begin{equation}
\label{DecompositionEq}
\Tm{\Z}\oplus (\Spec E)(1)\oplus\Tm{\Z}(2)
\end{equation}
for every $k$-\'etale algebra~$E$, such that $\Pic{\bar{S}}$
is a direct $G_{k}$-summand of~$\CH_{0}(\bar{E})$. Note that the
summand~$\Tm{\Z}$ corresponds to the idempotent
$\eta\times [S]\in\CH_{2}(S\times_{k}S)$, and the summand~$\Tm{\Z}(2)$
to the idempotent $[S]\times\eta\in\CH_{2}(S\times_{k}S)$, where~$\eta\in\CH_{0}(S)$
has degree~$1$. In particular the ``middle part''~$(S,\rho)$ of the
motive of~$S$, where $\rho=\id_{S}-(\eta\times [S]+[S]\times\eta)$,
is a direct summand of~$(\Spec E)(1)$.

\smallbreak

In some cases the complement of the summand~$S$ in~(\ref{DecompositionEq})
can be computed explicit. Recall for this that a $G_{k}$-module~$C$ is called {\it coflabby}
if~$\HM^{1}(H,C)=0$ for all open subgroups $H\subseteq G_{k}$.
A {\it coflabby resolution} of~$\Pic\bar{S}$ is an exact sequence
\begin{equation}
\label{CoflasqueResEq}
\xymatrix{
0 \ar[r] & C \ar[r]^-{\iota} & P \ar[r]^-{f} & \Pic\bar{S} \ar[r] & 0
}
\end{equation}
with~$P$ a $G_{k}$-permutation module and~$C$ a
coflabby $G_{k}$-module. Such a resolution always exists,
and splits if and only if $\Pic\bar{S}$ is an invertible
$G_{k}$-module, \ie a direct summand of a $G_{k}$-permutation
module, see \eg~\cite[Lem.\ 1]{CTSa77a}. Hence the motive of~$S$
is zero dimensional in~$\Chow (k)$ if and only if one (and so
all) coflabby resolutions of~$\Pic\bar{S}$ split.

\smallbreak

If there is such a split coflabby resolution~(\ref{CoflasqueResEq})
with~$C$ also a permutation module then
$$
(S,\rho)\oplus (\Spec F)(1)\,\simeq\, (\Spec E)(1)
$$
in~$\Chow (k)$, where~$E$ and~$F$ are a $k$-\'etale algebras, such
that $P\simeq\CH_{0}(\bar{E})$ and $C\simeq\CH_{0}(\bar{F})$ as
$G_{k}$-modules. We leave the details to the reader.
\end{emptythm}

\begin{emptythm}
\label{DelPezzoSubSect}
{\it Del Pezzo surfaces of degree~$5$.}
Let~$S$ be a Del Pezzo surface of degree~$d=5$.
Recall, see Manin's book~\cite[Sect.\ 25]{Cubic}, that~$\Pic\bar{S}=\CH_{1}(\bar{S})$
is a free $\Z$-module which has a basis $\ell_{0},\ell_{1},\ell_{2},\ell_{3},\ell_{4}$,
such that $(\ell_{i},\ell_{j})=0$ for $i\not= j$, $(\ell_{0},\ell_{0})=1$,
and $(\ell_{i},\ell_{i})=-1$ for $1\leq i\leq 4$, where~$(-,-)$ is the intersection pairing.
The class of the canonical bundle is then
$$
\varpi_{S}\, =\, -3\ell_{0}+\sum\limits_{i=1}^{4}\ell_{i}\in\Pic\bar{S}\, ,
$$
and
$$
\big\{\,\ell\in\Pic\bar{S}\, |\, (\varpi_{S},\ell)=0\, ,\; (\ell,\ell)=-2\,\big\}
$$
is a root system of type~${\rm A}_{4}$. A set of simple roots is given by
$s_{1}=\ell_{1}-\ell_{2}$, $s_{2}=\ell_{2}-\ell_{3}$, $s_{3}=\ell_{3}-\ell_{4}$,
and $s_{4}=\ell_{0}-\ell_{1}-\ell_{2}-\ell_{3}$.

\smallbreak

Identifying $\Z^{1+4}\xrightarrow{\simeq}\Pic\bar{S}$,
$(a,b_{1},b_{2},b_{3},b_{4})\longmapsto
a\ell_{0}-\sum\limits_{i=1}^{4}b_{i}\ell_{i}$,
the to the simple roots~$s_{i}$ corresponding
reflections~$\sigma_{i}$ act on~$\Pic\bar{S}$ as follows:
The map~$\sigma_{i}$ interchanges the coordinates~$b_{i}$ and~$b_{i+1}$
and keeps the other coordinates for $i=1,2,3$, and
$$
\sigma_{4}(a,b_{1},b_{2},b_{3},b_{4})=
(2a-b_{1}-b_{2}-b_{3},a-b_{2}-b_{3},a-b_{1}-b_{3},a-b_{1}-b_{2},b_{4})\, .
$$

\smallbreak

As the action of the Galois group~$G_{k}$ on~$\Pic\bar{S}$ fixes~$\varpi_{S}$
and preserves the intersection pairing this action factors through the above
described action of the Weyl group~${\rm W}({\rm A}_{4})$ of type~${\rm A}_{4}$,
see~\cite[Thm.\ 23.9]{Cubic}.

\bigbreak

Using these data we construct now a coflabby resolution for the $G_{k}$-module~$\Pic\bar{S}$.
Consider the elements $h_{i}:=\ell_{0}-\ell_{i}\in\Pic\bar{S}$,
$i=1,2,3,4$, and $h_{5}=2\ell_{0}-\ell_{1}-\ell_{2}-\ell_{3}-\ell_{4}$. Together with~$\varpi_{S}$
these elements generate the abelian group~$\Pic\bar{S}$. The reflection~$\sigma_{i}$
interchanges~$h_{i}$ and~$h_{i+1}$ and leaves the other $h_{j}$'s fixed for
$i=1,2,3,4$. Let
$$
P\, :=\;\Z\cdot e\,\oplus\,\bigoplus\limits_{i=1}^{5}\,\Z\cdot e_{i}
$$
a free $\Z$-module of rank~$6$. We let~${\rm W}({\rm A}_{4})$ act on this module as follows:
The simple reflection~$\sigma_{i}$ interchanges~$e_{i}$ and~$e_{i+1}$ and
leaves the other~$e_{j}$'s and~$e$ fixed. Hence it is a~${\rm W}({\rm A}_{4})$-permutation
module, and therefore also a $G_{k}$-permutation module. The map
$$
f\, :\; P\,\too\,\Pic\bar{S}\, ,
$$
which maps~$e_{i}$ to~$h_{i}$ and~$e$ to~$\varpi_{S}=
-3\ell_{0}+\sum\limits_{i=1}^{4}\ell_{i}$ is surjective, and ${\rm W}({\rm A}_{4})$-
and so also $G_{k}$-linear. Its kernel is generated by the ${\rm W}({\rm A}_{4})$- and so also
$G_{k}$-invariant element $x:=(e_{1}+e_{2}+e_{3}+e_{4}+e_{5})+2e$ and so is isomorphic to
the trivial $G_{k}$-module~$\Z$. Hence the exact sequence
\begin{equation}
\label{5resEq}
\xymatrix{
0 \ar[r] & {\Z} \ar[r]^-{\iota} &
P={\Z\cdot e\,\oplus\,\bigoplus\limits_{i=1}^{5}\,\Z\cdot e_{i}}
\ar[r]^-{f} & \Pic\bar{S} \ar[r] & 0\, ,
}
\end{equation}
where~$\iota$ maps~$1$ to~$x$, is a coflabby resolution of the
$G_{k}$-module~$\Pic\bar{S}$. It is split, as for instance $e\mapsto -2$
and $e_{i}\mapsto 1$ for~$i=1,2,3,4,5$ splits~$\iota$.

\medbreak

Since a Del Pezzo surface of degree~$5$ always has a rational point
(in fact is rational) by a classical result of Enriques, see for
instance Skorobogatov~\cite[Sect.\ 3.1]{TRP} for a ``modern'' proof,
the considerations above imply the following computation of motives.
\medbreak

\noindent
{\bfseries Theorem.}
{\it
Let~$S$ be a Del-Pezzo surface of degree~$5$ over~$k$. Then
we have an isomorphism
\begin{equation}
\label{DelPezzo5Eq}
\Tm{\Z}(1)\oplus S\,\simeq\,\Tm{\Z}\oplus\Tm{\Z}(1)\oplus (\Spec E)(1)
\oplus\Tm{\Z}(2)
\end{equation}
in~$\Chow (k)$, where~$E$ is a $k$-\'etale algebra of degree~$5$.
}

\medbreak

\noindent
{\bfseries Remarks.}
\begin{itemize}
\item[(i)]
A small alteration of the construction above gives a coflabby resolution
for the Picard group~$\Pic\bar{S}$ where~$S$ is a Del Pezzo surface of
degree~$6$ (we have essentially only to remove~$\ell_{4}$).

\smallbreak

\item[(ii)]
The summand~$\Tm{\Z}(1)$ on both sides of the isomorphism~(\ref{DelPezzo5Eq})
does not cancel if~$\Pic\bar{S}$ is not a $G_{k}$-permutation module.
\end{itemize}
\end{emptythm}

\begin{emptythm}
\label{torsionCH0SubSect}
{\it Torsion in the group of zero cycles of a geometrically rational surface.}
Let~$S$ be a geometrically rational $k$-surface.
If~$S$ is zero dimensional in~$\Chow (k)$ then~$S$ is strongly
$\Z$-trivial, see the example in~\ref{Z-trivialSubSect}.

\smallbreak

The following
theorem of Colliot-Th\'el\`ene~\cite{BriefCT13}, see the Appendix,
implies that the converse is also true.

\medbreak

\noindent
{\bfseries Theorem (Colliot-Th\'el\`ene).}
{\it
Let~$K$ be a field with separable closure~$K_{s}$
and~$X$ a smooth projective and geometrically integral
$K$-scheme with $K$-rational point, such that
\begin{itemize}
\item[(i)]
$\Pic X_{K_{s}}$ is a free and finitely generated $\Z$-module, and

\smallbreak

\item[(ii)]
$\deg :\CH_{0}(X_{K(X)})\too\Z$ is an isomorphism.
\end{itemize}
Then $\Pic X_{K_{s}}$ is an invertible $G_{K}$-module,
where~$G_{K}=\Gal (K_{s}/K)$ is the absolute Galois
group of~$K$.
}

\medbreak

\noindent
(Note that if~$K$ is a perfect field the assumption~$X(K)\not=\emptyset$
can be replaced by~$X$ has a $0$-cycle of degree~$1$ as under the assumption~$K$
is perfect the existence of such a cycle assures already that~$X$ possesses
a universal torsor, see~\cite[(2.0.2), and Props. 2.2.2~(b) and 2.2.5]{CTSa87}.)

\medbreak

Let~$S$ be a geometrically rational $k$-surface with function
field~$K=k(S)$, which is strongly $\Z$-trivial. Then also~$S_{K}$
is strongly $\Z$-trivial and so by the theorem of Colliot-Th\'el\`ene
the $G_{K}$-module $\Pic S_{K_{s}}$ is invertible. But $\Gal (K_{s}/\bar{k}(S))$
acts trivially on $\Pic S_{K_{s}}$ since the restriction maps
$\Pic\bar{S}\too\Pic S_{\bar{k}(S)}\too\Pic S_{K_{s}}$ are both
isomorphisms. Identifying~$G_{k}$ with~$\Gal (\bar{k}(S)/K)$ it
follows that also the $G_{k}$-module $\Pic\bar{S}$ is invertible.
Therefore:
\end{emptythm}

\begin{emptythm}
\label{mainThm}
{\bfseries Theorem.}
{\it
Let~$S$ be a geometrically rational surface over~$k$.
Then the following statements are equivalent:
\begin{itemize}
\item[(i)]
The degree map $\CH_{0}(S_{k(S)})\too\Z$ is an isomorphism
and~$S$ has a $0$-cycle of degree~$1$.

\smallbreak

\item[(ii)]
The surface~$S$ is strongly $\Z$-trivial.

\smallbreak

\item[(iii)]
The $G_{k}$-module~$\Pic\bar{S}$ is invertible
and~$S$ has a $0$-cycle of degree~$1$.

\smallbreak

\item[(iv)]
The motive of~$S$ in~$\Chow (k)$ is zero dimensional.
\end{itemize}
}
\end{emptythm}

\begin{emptythm}
\label{CH0-SubSect}
{\it Strong $\Z$-triviality of $T^{S}$-torsors over~$S$.}
The following is a corollary of the arguments in the letter~\cite{BriefCT13} by
Colliot-Th\'el\`ene.

\smallbreak
 
Let~$S$ be a geometrically rational surface over~$k$, whose motive
in~$\Chow (k)$ is $0$-dimensional, and~$T^{S}$ the torus with character
group $\Pic\bar{S}$.

\smallbreak

Let~$Y$ be a $T^{S}$-torsor over~$S$, and~$Y^{c}$ and $(T^{S})^{c}$ smooth
compactifications of~$Y$ and~$T^{S}$, respectively. Such compactifications exists,
see Colliot-Th\'el\`ene, Harari, and Skorobogatov~\cite{CTHaSk05},
and Colliot-Th\'el\`ene and Sansuc~\cite[Rem.\ 2.1.4]{CTSa87}. We aim to
show that under the assumption that the Chow motive of~$S$ is
$0$-dimensional the $k$-scheme~$Y^{c}$ is strongly $\Z$-trivial, \ie
$$
\deg\, :\;\CH_{0}(Y^{c}_{K})\,\too\,\Z
$$
is an isomorphism for all field
extensions $K\supseteq k$ (\cf~\ref{Z-trivialSubSect}).

\smallbreak

For this we observe first that since the motive of~$S$ is $0$-dimensional
in~$\Chow (k)$ it is a direct summand of $\Tm{\Z}\oplus (\Spec E)(1)\oplus\Tm{\Z}(2)$
for some $k$-\'etale algebra~$E$, see~\ref{ComputationSubSect}. Hence~$S_{K}$
is a direct summand of $\Tm{\Z}\oplus (\Spec E_{K})(1)\oplus\Tm{\Z}(2)$ in~$\Chow (K)$
and so it is also $0$-dimensional for all fields~$K\supseteq k$.

\smallbreak

By Theorem~\ref{mainThm} above we know that then $\Pic\bar{S}$ is
an invertible $G_{k}$-module and therefore~$T^{S}$ is a direct
summand of a quasi trivial torus, which in turn by Hilbert~90
implies that $\HM^{1}(K,T^{S}_{K})=0$ for all $K\supseteq k$. Hence
every $T^{S}$-torsor over a field extensions of~$k$ is trivial and
so in particular $k(S)\times_{S}Y$ is trivial. It follows that~$Y$
is birational to $S\times_{k}T^{S}$. Then also $Y^{c}_{K}$ is birational
isomorphic to $(S\times_{k}(T^{S})^{c})_{K}$ for all~$K\supseteq k$,
and so we have by the lemma in~\ref{ChowGroupSubSect} that 
$\deg :\CH_{0}(Y^{c}_{K})\too\Z$ is an isomorphism if and only
if $\deg :\CH_{0}(S_{K}\times_{K}(T^{S})^{c}_{K})\too\Z$ is an isomorphism.
But since the Chow motive of~$S_{K}$ is a direct summand of
$\Tm{\Z}\oplus (\Spec E_{K})(1)\oplus\Tm{\Z}(2)$ the latter group
is isomorphic to~$\CH_{0}((T^{S})^{c}_{K})$ via the push-forward along
the projection $S_{K}\times_{K}(T^{S})^{c}_{K}\too (T^{S})^{c}_{K}$ for
all field extensions~$K$ of~$k$. Hence it is enough to show that
the degree map
$$
\deg\, :\;\CH_{0}((T^{S})^{c}_{K})\,\too\,\Z
$$
is an isomorphism for all field extensions $K\supseteq k$.

\smallbreak

Let~$K$ be such an extension with separable closure~$K_{s}$ and
absolute Galois group $G_{K}=\Gal (K_{s}/K)$.
The projections from $S_{K}\times_{K} (T^{S})^{c}_{K}$ to~$S_{K}$
and to~$(T^{S})^{c}_{K}$ induce a $G_{K}$-linear map
$$
\Pic S_{K_{s}}\oplus\Pic (T^{S})^{c}_{K_{s}}\,\too\,
\Pic\big(\, S_{K_{s}}\times_{K_{s}} (T^{S})^{c}_{K_{s}}\,\big)\, ,
$$
which is an isomorphism by Colliot-Th\'el\`ene and Sansuc~\cite[Lem.\ 11]{CTSa77a}
since~$S_{K_{s}}$ is rational as shown by Coombes~\cite{Co88}.
Therefore~$\Pic (T^{S})^{c}_{K_{s}}$ is a
direct $G_{K}$-summand of $\Pic (S_{K_{s}}\times_{K_{s}} (T^{S})^{c}_{K_{s}})$.

\smallbreak

On the other hand by~\cite[Prop.\ 2.A.1]{CTSa87} there are
$G_{K}$-permutation modules~$P$ and~$Q$, such that we have an isomorphism
of $G_{K}$-modules
\begin{equation}
\label{G-IsomEq}
\Pic (S_{K_{s}}\times_{K_{s}} (T^{S})^{c}_{K_{s}})\oplus P\,\simeq\,
\Pic Y^{c}_{K_{s}}\oplus Q
\end{equation}
since~$Y^{c}_{K}$ and $S_{K}\times_{K} (T^{S})^{c}_{K}$ are $K$-birational
to each other.

\smallbreak

This is in particular true for an universal $T^{S}$-torsor~$Z$ over~$S$. For such
a $T^{S}$-torsor the $G_{K}$-module~$\Pic Z^{c}_{K_{s}}$ is invertible
by~\cite[Thm.\ 2.1.2]{CTSa87} and therefore by the isomorphism~(\ref{G-IsomEq})
for~$Y=Z$ also $\Pic (T^{S})^{c}_{K_{s}}$ is an
invertible $G_{K}$-module. But this implies by Colliot-Th\'el\`ene and
Sansuc~\cite[Prop.\ 6]{CTSa77a} (\cf the remark below) that~$T^{S}_{K}$ and so also its
compactification~$(T^{S})^{c}_{K}$ is stably rational, and so
the degree map $\deg :\CH_{0}((T^{S})^{c}_{K})\too\Z$ is an isomorphism
as claimed. We have proven one direction of the following result.
\end{emptythm}

\begin{emptythm}
\label{universaltorsorCor}
{\bfseries Corollary.}
{\it
Let~$S$ be a geometrically rational $k$-surface. Denote
by~$T^{S}$ the $k$-torus with character group~$\Pic\bar{S}$. Then the
following is equivalent:
\begin{itemize}
\item[(i)]
The motive of~$S$ in~$\Chow (k)$ is $0$-dimensional.

\smallbreak

\item[(ii)]
Every compactification of every $T^{S}$-torsor over~$S$
is strongly $\Z$-trivial.
\end{itemize}
}

\begin{proof}
We are left to show (ii)~$\Longrightarrow$~(i).
Assuming~(ii) the scheme $S\times_{k}(T^{S})^{c}$ is strongly
$\Z$-trivial, where~$(T^{S})^{c}$ is a compactification of~$T^{S}$,
and so the degree map
$\CH_{0}\big((S\times_{k}(T^{S})^{c})_{K}\big)\too\Z$
is an isomorphism for all field extensions $K\supseteq k$.
Since this degree map factors through
$\deg :\CH_{0}(S_{K})\too\Z$ also the latter is surjective.
On the other hand~$T^{S}_{K}$ and so also $(T^{S})^{c}_{K}$ has
a $K$-rational point and therefore~$\CH_{0}(S_{K})$ is
a direct summand of $\CH_{0}\big((S\times_{k}(T^{S})^{c})_{K}\big)$
and consequently torsion free. It follows that $\deg :\CH_{0}(S_{K})\too\Z$
is an isomorphism as the kernel of this map is the torsion subgroup
of~$\CH_{0}(S_{K})$ since~$S$ is geometrically rational. Hence~$S$ is
strongly $\Z$-trivial and so the motive of~$S$ in~$\Chow (k)$ is
$0$-dimensional by Theorem~\ref{mainThm} above. We are done.
\end{proof}

\smallbreak

\noindent
{\bfseries Remark.}
The assertion of~\cite[Prop.\ 6]{CTSa77a} we use assumes that the base field has
characteristic~$0$ to assure the existence of a smooth compactification
for every torus over this field. Only later it has been shown by
Brylinski and K\"unnemann that such compactifications
exist also in positive characteristic, \cf~\cite{CTHaSk05}.
\end{emptythm}

\goodbreak
\appendix
\section{A letter of Colliot-Th\'el\`ene to the author}
\bigbreak

\noindent
The following is a reproduction of the letter~\cite{BriefCT13} by
Colliot-Th\'el\`ene except that (i)~the list of references has been
deleted, (ii)~all citations refer now to the bibliography of the article,
and~(iii) a remark at the end has been added.

\bigbreak

Toronto,  den 8. Mai 2013

\medbreak

Brief von J.-L. Colliot-Th\'el\`ene an Stefan Gille

\bigbreak

Hier ist ein allgemeiner Satz.

\medbreak

\noindent
{\bf Satz.}
{\it 
Sei $F$ ein K\"orper, $X$ eine glatte, projektive, geometrisch irreduzible Variet\"at 
mit einem rationalen Punkt. Sei $G=\Gal(F_{s}/F)$ die absolute Galoisgruppe von $F$.

Nehmen wir an :

(a)  F\"ur jede Erweiterung von K\"orpern $L/F$ 
ist die Gradabbildung
$$
\CH_{0}(X_{L})\too \Z
$$
bijektiv. 

(b)
$M={\rm Pic}(X\times_{F}F_{s})$  ist ein Gitter.

Dann gibt es ein $G$-Gitter $N$, einen Permutationsmodul $P$ und einen
Isomorphismus
$$
M \oplus N \simeq P
$$
von Galois Gittern.}

\medskip

Beweis -- auch im Falle ${\rm Char}(F)>0$.  

Sei $p_{0} \in X(F)$, $T$ der $F$-Torus mit Charaktergruppe $\hat T = M$, und
$Y \too X$ der universelle Torsor \"uber $X$ mit trivaler Faser im Punkt $p_{0}$.
Dies ist ein Torsor unter der Gruppe $T$. Sei $L/F$ eine Erweiterung von K\"orpern.
So ein Torsor definiert Abbildungen
$$
X(L) \too\HM^{1}(L,T)
$$
und allgemeiner
$$
\CH_{0}(X_{L}) \too\HM^{1}(L,T)
$$
(letztere eine Gruppenabbildung)
die $p_{0}$ nach $0$ schicken.
(Siehe~\cite[Prop.\ 12 Seite 198]{CTSa77a}).

Sei $\eta$ der generische Punkt von $X$ und $L=F(X)$ der Funktionenk\"orper von $X$.
Unter der Annahme $Grad : \CH_{0}(X_{L}) \simeq  \Z$ folgt, da{\ss} das Bild von $\eta - p_{0}$ unter
$$\
CH_{0}(X_{L}) \too\HM^{1}(L,T)
$$ null ist, also ist   das Bild von~$\eta$ unter
$$
X(L)\too\HM^{1}(L,T)
$$
auch null. Das besagt aber, da{\ss} der Torsor $Y\too X$
unter $T$ einen rationalen Schnitt hat, und somit ist die $F$-Variet\"at $Y$
$F$-birational zum Produkt $X \times_{F}T$.

Sei $T^{c}/F$ eine \"aquivariante glatte Kompaktifizierung vom $F$-torus $T$
(\cite{CTHaSk05}).

Dann ist $Y^{c}=Y\times^{T}T^{c}$ (contracted product) eine vollst\"andige
glatte Kompaktifizierung von $Y$.

Die glatten projektiven Variet\"aten $Y^{c}$ und $X \times_{F}T^{c}$ sind $F$-birational zueinander.

In~\cite[Th\'eor\`eme 1]{CTSa77b} und~\cite[Thm. 2.12 Seite 411]{CTSa87}  wird gezeigt : 
Da $Y\too X$  ein universeller Torsor ist,  folgt
$F_{s} \simeq F_{s}[Y]^{\times}$ und $\Pic (Y\times_{F}F_{s})=0$.
Daraus folgt leicht : $\Pic (Y^{c}_{F_{s}})$ ist ein Gitter und als $G$-Modul
ein Permutationsmodul $P_{0}$.

Weil $T^{c}$ geometrisch rational ist, ist die nat\"urliche Abbildung
$$
\Pic (X_{F_{s}})\oplus\Pic (T^{c}_{F_{s}}) \too\Pic((X \times_{F}T^{c})_{F_{s}})
$$
ein Isomorphismus (\cite[Lemme 11 Seite 188]{CTSa77a}).

\medbreak

In \cite[Prop. 2.A.1 Seite 461]{CTSa87} 
findet man einen Beweis von Moret-Bailly f\"ur das folgende Lemma.
  
\medbreak

\noindent
{\bfseries Lemma.}
{\it
Seien $V$ und $W$ zwei glatte vollst\"andige geometrisch irreduzible Variet\"aten
die $F$-birational zueinander sind. Dann gibt es Permutationsmodule $P_{1}$ und $P_{2}$,
und einen Isomorphismus von Galoismodulen
$$
P_{1} \oplus \Pic (V_{F_{s}}) \simeq P_{2} \oplus \Pic (W_{F_{s}})\, .
$$
}
 
\bigbreak

Wenn man alles zusammenfasst, erh\"alt man einen Isomorphismus von Galoismodulen
$$
P_{1} \oplus  P_{0}  \simeq P_{2} \oplus \Pic (X_{F_{s}})  \oplus\Pic (T^{c}_{F_{s}})\, ,
$$
also ist  ${\rm Pic}(X_{F_{s}})$ invertierbar.
 
Was zu beweisen war.

\bigbreak

\noindent
{\bfseries Bemerkungen.}

\medbreak

1. Wenigstens wenn ${\rm Char}(F)=0$ ist,  sollte man zeigen k\"onnen,
da{\ss}  Annahme (b) aus Annahme (a) folgt. Man kann schon sehen,
da{\ss}  $\Pic (X\times_{F}F_{s})$     von endlichem Typ \"uber $\Z$ ist.
 
(Hinzugef\"ugt, 10. April 2014) Sei ${\rm Char}(F)=0$, und sei $n>0$ eine ganze Zahl.
Aus  Annahme (a)   und der  Kummerschen Sequenz folgt :
$$
0= \HM^{1}_{nr}(F_{s}(X)/F_{s},\mu_{n})=\HM^{1}_{\et}(X_{s},\mu_{n})=\Pic (X_{F_{s}})[n]\, .
$$
Wegen der bekannten Struktur der Picard Gruppe folgt, da{\ss} die Picard Variet\"at von 
$X$ null ist, und da{\ss} die N\'eron-Severi Gruppe torsionsfrei ist, also ist 
$\Pic (X\times_{F}F_{s})$  ein Gitter.
  
 \medbreak
 
2.  Es ist nicht klar, ob die Arbeit von Merkurjev~\cite{Me08b}
einen alternativen Beweis des Satzes geben k\"onnte. Aus dieser Arbeit
kann man wenigstens folgern, da{\ss} $\Pic (X\times_{F}F_{s})$  co-welk ist.

\bigbreak\bigbreak

\bibliographystyle{amsalpha}

\begin{thebibliography}{99}

\smallbreak
\bibitem{StAlg}
A.\ A.\ Albert, \textsl{Structure of Algebras}, American Mathematical
Society Colloquium Publications, Vol.\ {\bfseries 24}, American Mathematical
Society, New-York, 1939.

\smallbreak
\bibitem{Bl81}
S.\ Bloch, \textsl{On the Chow groups of certain rational surfaces},
Ann.\ Sci.\ \'Ecole Norm.\ Sup.\ (4) {\bfseries 14} (1981), 41--59.

\smallbreak
\bibitem{ChMe06}
V. Chernousov, A.\ S.\ Merkurjev, \textsl{Motivic decomposition of projective
homogeneous varieties and the Krull-Schmidt theorem}, Transform.\ Groups
{\bfseries 11} (2006), 371--386.

\smallbreak
\bibitem{CT83}
J.-L.\ Colliot-Th\'el\`ene, \textsl{Hilbert's Theorem 90
for $\QK_{2}$, with application to the Chow groups of
rational surfaces}, Invent.\ Math.\ {\bfseries 71} (1983), 1--20.

\smallbreak
\bibitem{BriefCT13}
J.-L.\ Colliot-Th\'el\`ene, Letter to the author, dated May 8, 2013
(in German).

\smallbreak
\bibitem{CTHaSk05}
J.-L.\ Colliot-Th\'el\`ene, D.\ Harari, A.\ N.\ Skorobogatov, \textsl{Compactification
\'equivariante d'un tore (d'apr\`es Brylinski et K\"unnemann)}, Expo.\ Math.\
{\bfseries 23} (2005), 161--170.

\smallbreak
\bibitem{CTSa77a}
J.-L.\ Colliot-Th\'el\`ene, J.-J.\ Sansuc, \textsl{La $R$-\'equivalence sur les
tores}, Ann.\ Sci.\ \'Ecole Norm Sup.\ (4) {\bfseries 10} (1977), 175--229.

\smallbreak
\bibitem{CTSa77b}
J.-L.\ Colliot-Th\'el\`ene, J.-J.\ Sansuc, \textsl{Vari\'et\'es de premi\`ere descente
attach\'ees aux vari\'et\'es rationnelles}, C.\ R.\ Acad.\ Sci.\ Paris S\'er.\ A-B
{\bfseries 284} (1977), A967--A970.

\smallbreak
\bibitem{CTSa87}
J.-L.\ Colliot-Th\'el\`ene, J.-J.\ Sansuc, \textsl{La descente sur les vari\'et\'es
rationnelles. II}, Duke Math.\ J.\ {\bfseries 54} (1987), 375--492.

\smallbreak
\bibitem{Co88}
K.\ R.\ Coombes, \textsl{Every rational surface is separably split}, Comment.\
Math.\ Helv.\ {\bfseries 63} (1988), 305--311.

\smallbreak
\bibitem{IntTheory}
W.\ Fulton, \textsl{Intersection theory}, Springer-Verlag, Berlin, 1984.

\smallbreak
\bibitem{Gi10}
S.\ Gille, \textsl{The Rost nilpotence theorem for
geometrically rational surfaces}, Invent.\ Math.\
{\bfseries 181} (2010), 1--19.

\smallbreak
\bibitem{Gi14}
S.\ Gille, \textsl{On Chow motives of surfaces}, J.\ reine angew.\
Math.\ {\bfseries 686} (2014), 149--166.

\smallbreak
\bibitem{KaMe13}
N.\ Karpenko, A.\ S.\ Merkurjev, \textsl{On standard norm varieties}, Ann.\
Sci.\ \'Ecole Norm.\ Sup.\ (4) {\bfseries 46} (2013), 175--214.

\smallbreak
\bibitem{Ma68}
Y.\ Manin, \textsl{Correspondences, motifs and monoidal transformations},
(Russian), Mat.\ Sb.\ (N.S.) {\bfseries 77} (1968), 475--507.

\smallbreak
\bibitem{Cubic}
Y.\ Manin, \textsl{Cubic forms: algebra, geometry, arithmetic},
North-Holland Publishing 1974.

\smallbreak
\bibitem{Me08a}
A.\ S.\ Merkurjev, \textsl{$R$-equivalence on three-dimensional tori and
zero-cycles}, Algebra $\&$ Number Theory {\bfseries 2} (2008), 69--89.

\smallbreak
\bibitem{Me08b}
A.\ S.\ Merkurjev, \textsl{Unramified elements in cycle modules}, J.\ Lond.\
Math.\ Soc.\ (2) {\bfseries 78} (2008), 51--64.

\smallbreak
\bibitem{Ro96}
M.\ Rost, \textsl{Chow groups with coefficients}, Doc.\ Math.\
{\bfseries 1} (1996), 319--393.

\smallbreak
\bibitem{TRP}
A.\ N.\ Skorobogatov, \textsl{Torsors and rational points}, Cambridge
Univ.\ Press 2001.


\end{thebibliography}

\end{document}